\documentclass[11pt]{article}
\usepackage{amsmath}
\usepackage{amssymb}
\usepackage{amsthm}
\usepackage{cases}
\usepackage{url}

\newcommand\nonu{\nonumber}
\newcommand\dstyle\displaystyle
\newcommand\sa{\smallskipamount}
\newcommand\ma{\medskipamount}

\newcommand\sLP{\\[\sa]}

\newcommand\sPP{\\[\sa]\indent}
\newcommand\mPP{\\[\ma]\indent}
\newcommand\CC{\mathbb{C}}
\newcommand\RR{\mathbb{R}}
\newcommand\FSP{{\cal P}}
\newcommand\al\alpha
\newcommand\be\beta
\newcommand\de\delta
\newcommand\tha\theta
\newcommand\la\lambda
\newcommand\Ga{\Gamma}
\newcommand\De{\Delta}
\newcommand\half{\frac12}
\newcommand\thalf{\tfrac12}
\newcommand\iy\infty
\newcommand\wt{\widetilde}
\newcommand\union{\cup}
\newcommand\lan{\langle}
\newcommand\ran{\rangle}
\newcommand\const{{\rm const.}\,}
\newcommand{\hyp}[5]{\,\mbox{}_{#1}F_{#2}\!\left(
  \genfrac{}{}{0pt}{}{#3}{#4};#5\right)}
\newcommand{\qhyp}[5]{\,\mbox{}_{#1}\phi_{#2}\!\left(
  \genfrac{}{}{0pt}{}{#3}{#4};#5\right)}
\newcommand\LHS{left-hand side}
\newcommand\RHS{right-hand side}

\makeatletter
\newcommand{\raisemath}[1]{\mathpalette{\raisem@th{#1}}}
\newcommand{\raisem@th}[3]{\raisebox{#1}{$#2#3$}}
\makeatother

\numberwithin{equation}{section}
\newtheorem{theorem}{Theorem}[section]
\newtheorem{example}[theorem]{Example}
\newcommand\Proof{\noindent{\bf Proof}\quad}

\begin{document}
\title{Orthogonal polynomials, a short introduction}
\author{Tom H. Koornwinder}
\date{}
\maketitle
\begin{abstract}
This paper is a short introduction to orthogonal polynomials, both the general
theory and some special classes. It ends with some remarks about the
usage of computer algebra for this theory.
The paper has appeared as a chapter ``Orthogonal polynomials''
in the book
``Computer Algebra in Quantum Field Theory'', Springer Vienna, 2013. In the
present arXiv version some minor corrections were made.
\end{abstract}
%
%
\section{Introduction}
This paper is a short introduction to orthogonal polynomials, both the general
theory and some special classes. After the definition and first examples in
Section 2, important but mainly elementary aspects of the general theory
associated with the three-term recurrence relation are treated in Section 3.
Sections 4, 6 and 7 discuss special classes of orthogonal polynomials,
interrupted by Section 5 about Gauss quadrature.
Section  8 collects some more advanced results in the general
theory of orthogonal polynomials. Finally Section 9 discusses the role
of computer algebra in the theory of (special) orthogonal polynomials.

Everything treated here is well-known from the literature. I mention a few books
which can be recommended for more detailed study.
A great classical introduction to orthogonal polynomials, both the general theory
and the special polynomials, is Szeg\H o \cite{3}.
A very readable textbook, in particular for the general theory, is
Chihara \cite{1}. As a textbook emphasizing the special theory I recommend
Andrews, Askey \& Roy \cite{2}. Very good is also Ismail \cite{6}, but more
focusing on the $q$-case. Two recent compendia of formulas for special
orthogonal polynomials are Olver et al.\ \cite[Ch. 18]{5} and
Koekoek, Lesky \& Swarttouw \cite[Ch. 9 and 14]{4}.

\section{Definition of orthogonal polynomials and first examples}
Let $\FSP$ be the real vector space of all polynomials in one variable with
real coefficients.
Assume on $\FSP$ a (positive definite) inner product $\lan f,g\ran$ ($f,g\in\FSP$).
Orthogonalize the sequence pf monomials
$1,x,x^2,\ldots\;$ with respect to the inner product
(Gram--Schmidt). This results into the sequence
$p_0(x),p_1(x),p_2(x),\ldots\;$ of polynomials in $x$.
So $p_0(x)=1$ and, if $p_0(x),p_1(x),\ldots,p_{n-1}(x)$ are already produced
and mutually orthogonal, then
\[
p_n(x):=x^n-\sum_{k=0}^{n-1}\frac{\lan x^n,p_k\ran}{\lan p_k,p_k\ran}\,p_k(x).
\]
Indeed, $p_n(x)$ is a linear combination of $1,x,x^2,\ldots,x^n$, and
\begin{equation*}
\lan p_n,p_j\ran=\lan x^n,p_j\ran
-\sum_{k=0}^{n-1}\frac{\lan x^n,p_k\ran}{\lan p_k,p_k\ran}\,\lan p_k,p_j\ran
=\lan x^n,p_j\ran-\frac{\lan x^n,p_j\ran}{\lan p_j,p_j\ran}\,\lan p_j,p_j\ran=0
\quad(j=0,1,\ldots,n-1).
\end{equation*}

Throughout we will use the constants $h_n$ and $k_n$ associated with the
orthogonal system:
\begin{equation}
\lan p_n,p_n\ran=h_n,\qquad
p_n(x)=k_n x^n+\mbox{polynomial of lower degree}.
\label{20}
\end{equation}
The $p_n$ are unique up to a nonzero constant real factor. We may take them,
for instance,
{\em orthonormal} ($h_n=1$; this determines $p_n$ uniquely if also $k_n>0$)
or {\em monic} ($k_n=1$).

In general we want
\[
\lan x\,f,g\ran=\lan f,x\,g\ran.
\]
This is true, for instance, if
\[
\lan f,g\ran:=\int_a^b f(x)\,g(x)\,w(x)\,dx\qquad{\rm or}\qquad
\lan f,g\ran:=\sum_{j=0}^\iy f(x_j)\,g(x_j)\,w_j
\]
for a {\em weight function} $w(x)\ge0$
or for {\em weights} $w_j>0$, respectively.
These are special cases of an inner product
\begin{equation}
\lan f,g\ran:=\int_\RR f(x)\,g(x)\,d\mu(x),
\label{9}
\end{equation}
where $\mu$ is a (positive) {\em Borel measure} on $\RR$,
namely the cases
$d\mu(x)=w(x)\,dx$ on an interval $I$,
and $\mu=\sum_{j=1}^\iy w_j\,\de_{x_j}$, respectively.

A Borel measure $\mu$ on $\RR$ such that $\mu(K)<\iy$ for $K$ compact,
can also be thought of as a \emph{Lebesgue--Stieltjes measure}
\cite[Theorem 1.16]{27}. This involves a
{\em non-decreasing function}
$\wt \mu$ on $\RR$ such that, for suitable functions $f$,
we have a {\em Riemann--Stieltjes integral}
\[
\int_\RR f(x)\,d\mu(x)=\int_\RR f(x)\,d\wt\mu(x)=
\lim_{M\to\iy}\int_{-M}^M f(x)\,d\wt\mu(x).
\]
The measure $\mu$ has in $x$ a {\em mass point} of {\em mass} $c>0$ if
the non-decreasing function $\wt\mu$ has a jump $c$ at $x$,
i.e., if
$\lim_{\de\downarrow0}\,\big(\wt\mu(x+\de)-\wt\mu(x-\de)\big)=c>0$.
The number of mass points is countable.
More generally, the {\em support} of the measure $\mu$ consists of all
$x\in\RR$ such that
$\wt\mu(x+\de)-\wt\mu(x-\de)>0$ for all $\de>0$.
This set ${\rm supp}(\mu)$ is always a closed subset of $\RR$.

In the most general case let $\mu$ be a (positive) Borel measure on $\RR$
(of infinite support, i.e., not $\mu=\sum_{j=1}^Nw_j\,\de_{x_j}$)
such that
$\dstyle\int_\RR |x^n|\,d\mu(x)<\iy$ for all $n=0,1,2,\ldots\;$.
A system $\{p_0,p_1,p_2,\ldots\}$ obtained by orthogonalization of
$\{1,x,x^2,\ldots\}$ with respect to the inner product \eqref{9}
is called a
system of {\em orthogonal polynomials} with respect to the
orthogonality measure $\mu$.

Here follow some first examples of explicit orthogonal polynomials.
\begin{itemize}
\item
{\em Legendre polynomials} $P_n(x)$, orthogonal on $[-1,1]$ with respect
to the weight function 1. Normalized by $P_n(1)=1$.
\item
{\em Hermite polynomials} $H_n(x)$, orthogonal on $(-\iy,\iy)$ with
respect to the weight function $e^{-x^2}$.
Normalized by $k_n=2^n$.
\item
{\em Charlier polynomials} $c_n(x,a)$, orthogonal on the points $x=0,1,2,\ldots$ with respect to the weights $a^x/x!$ ($a>0$).
Normalized by $c_n(0;a)=1$.
\end{itemize}
The $h_n$ (see \eqref{20}) can be computed for these examples:
\begin{align*}
&\thalf\,\int_{-1}^1 P_m(x)\,P_n(x)\,dx=\frac1{2n+1}\,\de_{m,n}\,,\qquad
\pi^{-\half}\,\int_{-\iy}^\iy H_m(x)\,H_n(x)\,e^{-x^2}\,dx=2^nn!\,\de_{m,n}\,,
\\
&\qquad\qquad\qquad e^{-a}\,\sum_{x=0}^\iy c_m(x,a)\,c_n(x,a)\,\frac{a^x}{x!}
=a^{-n}n!\,\de_{m,n}\,.
\end{align*}

\section{Three-term recurrence relation and some consequences}
\subsection{Three-term recurrence relation}
The following theorem is fundamental for the general theory of orthogonal
polynomials.
\begin{theorem}
Orthogonal polynomials $p_n$ satisfy
\begin{equation}
\label{1}
\begin{split}
xp_n(x)&=a_n p_{n+1}(x)+b_n p_n(x)+c_n p_{n-1}(x)\quad(n>0),\\
xp_0(x)&=a_0 p_1(x)+b_0 p_0(x)
\end{split}
\end{equation}
with $a_n,b_n,c_n$ real constants and $a_nc_{n+1}>0$. Also
$\dstyle a_n=\frac{k_n}{k_{n+1}},\quad
\frac{c_{n+1}}{h_{n+1}}=\frac{a_n}{h_n}\,$.
\sPP
Moreover {\em (Favard theorem)},
if polynomials $p_n$ of degree $n$ $(n=0,1,2,\ldots)$ satisfy
\eqref{1}
with $a_n,b_n,c_n$ real constants and $a_nc_{n+1}>0$ then
there exists a {\rm(}positive{\rm)} measure $\mu$ on $\RR$ such that the polynomials
$p_n$ are orthogonal with respect to $\mu$.
\end{theorem}

The proof of the first part is easy. Indeed,
$xp_n(x)=\sum_{k=0}^{n+1}\al_k p_k(x)$,
and if $k\le n-2$ then
$\lan xp_n,p_k\ran=\lan p_n,xp_k\ran=0$, hence $\al_k=0$.
Furthermore, 
\[
c_{n+1}=\frac{\lan xp_{n+1},p_n\ran}{\lan p_n,p_n\ran}
=\frac{\lan xp_n,p_{n+1}\ran}{h_n}
=\frac{\lan xp_n,p_{n+1}\ran}{\lan p_{n+1},p_{n+1}\ran}\,
\frac{h_{n+1}}{h_n}=a_n\,\frac{h_{n+1}}{h_n}\,.
\]
Hence $a_nc_{n+1}=a_n^2\,h_{n+1}/h_n>0$.
Hence $c_{n+1}/h_{n+1}=a_n/h_n$.

The proof of the second part is much deeper (see Cihara \cite[Ch.~2]{1}).
\goodbreak
{\bf Remarks}
\begin{enumerate}
\item
For orthonormal polynomials the recurrence relation \eqref{1} becomes
\begin{equation}
\label{7}
\begin{split}
xp_n(x)&=a_n p_{n+1}(x)+b_n p_n(x)+a_{n-1} p_{n-1}(x)\quad(n>0),\\
xp_0(x)&=a_0 p_1(x)+b_0 p_0(x),
\end{split}
\end{equation}
and for monic orthogonal polynomials
\begin{equation}
\label{6}
\begin{split}
xp_n(x)&=p_{n+1}(x)+b_n p_n(x)+c_n p_{n-1}(x)\quad(n>0),\\
xp_0(x)&=p_1(x)+b_0 p_0(x),
\end{split}
\end{equation}
with $c_n=h_n/h_{n-1}>0$ in \eqref{6}.
If orthonormal polynomials $p_n$ satisfy \eqref{7} then the
corresponding monic polynomials $k_n^{-1}p_n$ satisfy \eqref{6} with
$c_n=a_{n-1}^2\,$.
\item
If the orthogonality measure is {\em even} ($\mu(E)=\mu(-E)$) then
$p_n(-x)=(-1)^n p_n(x)$,
hence $b_n=0$, so
$xp_n(x)=a_np_{n+1}(x)+c_np_{n-1}(x)$.
Examples of orthogonal polynomials with even orthogonality measure are the
Legendre and Hermite polynomials.
\item
The recurrence relation \eqref{1} determines the polynomials
$p_n$ uniquely (up to a constant factor because of the choice of the constant
$p_0$).
\item
The orthogonality measure $\mu$ for a system of orthogonal
polynomials may not be unique (up to a constant positive factor).
See Example \ref{21}.
\item
If $\mu$ is unique
then $\FSP$ is dense in $L^2(\mu)$. See Shohat \& Tamarkin
\cite[Theorem 2.14]{7}.
\item
If there is an orthogonality measure  $\mu$ with bounded support then $\mu$ is unique. See Chihara \cite[Ch.~2, Theorem~5.6]{1}.
\end{enumerate}

\subsection{Moments}
The {\em moment functional} $M\colon p\mapsto\lan p,1\ran\colon\FSP\to\RR$
associated with an orthogonality measure $\mu$ is already determined by the rule
$M(p_n)=\lan p_n,1\ran=0$ for $n>0$. Hence
$M$ is determined (up to a constant factor)
by the system of orthogonal polynomials $p_n$, independent of the choice
of the orthogonality measure, and hence $M$ is also determined by \eqref{1}.
The same is true for the inner product $\lan f,g\ran=\lan fg,1\ran$ on $\FSP$.

The moment functional $M$ is also determined by
the {\em moments} $\mu_n:=\lan x^n,1\ran$ ($n=0,1,2,\ldots$).
The condition $a_nc_{n+1}>0$ is equivalent to {\em positive definiteness} of the moments, stated as
\[
\De_n:=\det(\mu_{i+j})_{i,j=0}^n>0\quad(n=0,1,2,\ldots).
\]
For given moments $\mu_n$ and corresponding orthogonal polynomials $p_n$
a positive measure $\mu$ is an orthogonality measure for the $p_n$ iff
$\mu$ is a {\em solution of the {\rm(}Hamburger{\rm)} moment problem}
\begin{equation}
\int_\RR x^n\,d\mu(x)=\mu_n\qquad(n=0,1,2,\ldots).
\label{19}
\end{equation}
Uniqueness of the orthogonality measure is equivalent to uniqueness of the
moment problem.
\begin{example}[non-unique orthogonality measure]
\label{21}
\rm
The following goes back to Stieltjes \cite[\S56]{8}.
In the easily verified formula
\[
\int_{-\iy}^\iy e^{-u^2}(1+C\sin(2\pi u))\,du=\pi^{1/2}
\]
make a transformation of integration variable
 $u=\log x-\thalf(n+1)$ and take $-1<C<1$. Then
\begin{equation}
\pi^{-\half}e^{-\frac14}
\int_0^\iy x^n(1+C\sin(2\pi\log x))\,e^{-\log^2x}\,dx=e^{\frac14 n(n+2)}.
\label{22}
\end{equation}
Thus a one-parameter family of measures yields moments which
are independent of $C$. The corresponding orthogonal polynomials $p_n$
are a special case of the {\em Stieltjes-Wigert polynomials}
\cite[\S14.27]{4}:
$p_n(x)=S_n(q^\half x;q)$ with $q=e^{-\half}$, see
Christiansen \cite[p.223]{9}.

It is also elementary to show that
\begin{equation}
\left(\sum_{k=-\iy}^\iy e^{-\frac14 k^2}\right)^{-1}
\sum_{k=-\iy}^\iy e^{-\half kn} e^{-\frac14(k+1)^2}=e^{\frac14 n(n+2)},
\label{25}
\end{equation}
which means that the same moments, up to a constant factor, as in \eqref{22} are obtained with the measure $\sum_{k=-\iy}^\iy e^{-\frac14(k+1)^2} \de_{\exp(-\half k)}$.
\end{example}

\subsection{Christoffel--Darboux formula}
Let $\FSP_n$ be the space of polynomials of degree $\le n$. Let
$\{p_n\}$ be a system of orthogonal polynomials with respect to the measure
$\mu$. The
{\em Christoffel--Darboux kernel} is defined by
\begin{equation}
K_n(x,y):=\sum_{j=0}^n\frac{p_j(x)p_j(y)}{h_j}\,.
\label{31}
\end{equation}
Then
\[
(\Pi_nf)(x):=\int_\RR K_n(x,y)\,f(y)\,d\mu(y)
\]
defines an orthogonal projection $\Pi_n\colon \FSP\to\FSP_n\,$.
Indeed, if $f(y)=\sum_{k=0}^\iy \al_k p_k(y)$ (finite sum) then
\[
(\Pi_nf)(x)=\sum_{j=0}^n p_j(x)\sum_{k=0}^\iy \frac{\al_k}{h_j}
\int_\RR p_j(y)\,p_k(y)\,d\mu(y)
=\sum_{j=0}^n \al_j p_j(x).
\]
The {\em Christoffel--Darboux formula} for $K_n(x,y)$ given by \eqref{31}
is as follows.
\begin{numcases}
{\sum_{j=0}^n\frac{p_j(x)p_j(y)}{h_j}=}
\frac{k_n}{h_nk_{n+1}}\,\frac{p_{n+1}(x)p_n(y)-p_n(x)p_{n+1}(y)}{x-y}
&$(x\ne y)$,\label{5}\\
\frac{k_n}{h_nk_{n+1}}\,(p_{n+1}'(x)p_n(x)-p_n'(x)p_{n+1}(x))&$(x=y)$.
\label{10}
\end{numcases}
For the proof of \eqref{5} note that
\begin{align*}
xp_j(x)&=a_j p_{j+1}(x)+b_j p_j(x)+c_j p_{j-1}(x),\\
yp_j(y)&=a_j p_{j+1}(y)+b_j p_j(y)+c_j p_{j-1}(y).
\end{align*}
Hence
\begin{equation*}
\frac{(x-y)p_j(x)p_j(y)}{h_j}=
\frac{a_j}{h_j}(p_{j+1}(x)p_j(y)-p_j(x)p_{j+1}(y))
-\frac{c_j}{h_j}(p_j(x)p_{j-1}(y)-p_{j-1 }(x)p_j(y)).
\end{equation*}
Use that $c_j/h_j=a_{j-1}/h_{j-1}$. Sum from $j=0$ to $n$.
Use that $a_n=k_n/k_{n+1}$.
This yields \eqref{5}. For \eqref{10} let $y\to x$ in \eqref{5}.

\subsection{Zeros of orthogonal polynomials}
\begin{theorem}
Let $p_n$ be an orthogonal polynomial of degree $n$.
Let $\mu$ have support within the closure of the interval $(a,b)$.
Then $p_n$ has $n$ distinct zeros on $(a,b)$.
Furthermore, the zeros of $p_n$ and $p_{n+1}$ alternate.
\end{theorem}
\Proof
For the proof of the first part suppose $p_n$ has precisely
$k<n$ sign changes on $(a,b)$ at $x_1,x_2,\ldots,x_k$.
Hence, after possibly multiplying $p_n$ by $-1$, we have
$p_n(x)(x-x_1)\ldots(x-x_k)\ge0$ on $[a,b]$.
Hence $\int_a^b p_n(x)(x-x_1)\ldots(x-x_k)\,d\mu(x)>0$.
But by orthogonality we have\\
$\int_a^b p_n(x)(x-x_1)\ldots(x-x_k)\,d\mu(x)=0$.
Contradiction.\qed

For the proof of the second part use \eqref{10}: If $k_n,k_{n+1}>0$ then
\[
p_{n+1}'(x)p_n(x)-p_n'(x)p_{n+1}(x)
=\frac{h_nk_{n+1}}{k_n}\,\sum_{j=0}^n \frac{p_j(x)^2}{h_j}>0.
\]
Hence, if $y,z$ are two successive zeros of $p_{n+1}$ then
\[
p_{n+1}'(y)p_n(y)>0,\quad p_{n+1}'(z)p_n(z)>0.
\]
Since $p_{n+1}'(y)$ and $p_{n+1}'(z)$ have opposite signs,
$p_n(y)$ and $p_n(z)$ must have opposite signs.
Hence $p_n$ must have a zero in the interval $(y,z)$.\qed

\subsection{Kernel polynomials}
Recall the Christoffel-Darboux formula \eqref{5}.
Suppose the orthogonality measure $\mu$ has support within $(-\iy,b]$
and fix $y\ge b$. Then
\[
\int_{-\iy}^b K_n(x,y)\,x^k\,(y-x)\,d\mu(x)=y^k(y-y)=0
\qquad(k<n).
\]
Hence $x\mapsto q_n(x):=K_n(x,y)$ is an orthogonal polynomial
 of degree $n$ on $(-\iy,b]$ with respect
to the measure $(y-x)\,d\mu(x)$.
Hence
\begin{align*}
&q_n(x)-q_{n-1}(x)=\frac{p_n(y)}{h_n}\,p_n(x),\\
&p_n(y)p_{n+1}(x)-p_{n+1}(y)p_n(x)=\frac{h_nk_{n+1}}{k_n}\,(x-y)q_n(x).
\end{align*}
The orthogonal polynomials $q_n$ are called {\em kernel polynomials}.
Of course, they depend on the choice of $\mu$ and of $y$.

\section{Very classical orthogonal polynomials}
These are the Jacobi, Laguerre and Hermite polynomials.
They are usually called {\em classical orthogonal polynomials}, but I prefer
to call them {\em very classical} and to consider all polynomials in
the ($q$-)Askey scheme (see Sections 6 and 7) as classical.

We will need {\em hypergeometric series} \cite[Ch.~2]{2}:
\begin{equation}
\hyp rs{a_1,\ldots,a_r}{b_s,\ldots,b_s}z:=
\sum_{k=0}^\iy\frac{(a_1)_k\ldots (a_r)_k}{(b_1)_k\ldots(b_s)_k}\,
\frac{z^k}{k!}\,,
\label{34}
\end{equation}
where $(a)_k:=a(a+1)\ldots(a+k-1)$ for $k=1,2,\ldots\;$ and $(a)_0:=1$
is the {\em shifted factorial}. If one of the upper parameters in \eqref{34}
equals a
non-positive integer $-n$ then the series terminates after the term
with $k=n$.

\subsection{Jacobi polynomials}
{\em Jacobi polynomials} $P_n^{(\al,\be)}$
\cite[\S9.8]{4} are orthogonal on $(-1,1)$
with respect to the weight function
$w(x):=(1-x)^\al(1+x)^\be$ ($\al,\be>-1$) and they are normalized
by $P_n^{(\al,\be)}(1)=(\al+1)_n/n!\,$.
They can be expressed as terminating Gauss hypergeometric series:
\begin{align}
P_n^{(\al,\be)}(x)&=\frac{(\al+1)_n}{n!}\,
\hyp21{-n,n+\al+\be+1}{\al+1}{\thalf (1-x)}\nonu\\
&=\sum_{k=0}^n\frac{(n+\al+\be+1)_k\,\,(\al+k+1)_{n-k}}{k!\,(n-k)!}\,\frac{(x-1)^k}{2^k}\,.
\label{29}
\end{align}
They satisfy (because of the orthogonality property) the symmetry
$P_n^{(\al,\be)}(-x)=(-1)^n\,P_n^{(\be,\al)}(x)$.
Thus we conclude (much easier than by manipulation of the
hypergeometric series) that
\[
\hyp21{-n,n+\al+\be+1}{\al+1}z=
\frac{(-1)^n(\be+1)_n}{(\al+1)_n}\,
\hyp21{-n,n+\al+\be+1}{\be+1}{1-z}.
\]

For  $p_n(x):=P_n^{(\al,\be)}(x)$ there is the second order differential equation
\begin{equation}
(1-x^2)p_n''(x)+\big(\be-\al-(\al+\be+2)x\big)p_n'(x)=-n(n+\al+\be+1)\,p_n(x).
\label{30}
\end{equation}
This can be split up by the {\em shift operator relations}
\begin{align}
&\frac d{dx}\,P_n^{(\al,\be)}(x)=\thalf(n+\al+\be+1) P_{n-1}^{(\al+1,\be+1)}(x),
\label{11}\\
&(1-x^2)\,\frac d{dx}\,P_{n-1}^{(\al+1,\be+1)}(x)
+\big(\be-\al-(\al+\be+2)x\big)P_{n-1}^{(\al+1,\be+1)}(x)
\nonu\\
&\qquad
=\dstyle(1-x)^{-\al}(1+x)^{-\be}\,\frac d{dx}\left((1-x)^{\al+1}(1+x)^{\be+1}
P_{n-1}^{(\al+1,\be+1)}(x)\right)=-2n\,P_n^{(\al,\be)}(x).
\label{2}
\end{align}
Note that the operator $d/dx$ acting at the \LHS\ of \eqref{11} raises the
parameters and lowers the degree of the Jacobi polynomial,
while the operator acting at the \LHS\
of \eqref{2} lowers the parameters and raises the degree.
Iteration of \eqref{2} gives the {\em Rodrigues formula}
\[
 P_n^{(\al,\be)}(x)=\frac{(-1)^n}{2^n n!}\,(1-x)^{-\al}(1+x)^{-\be}
\left(\frac d{dx}\right)^n\left((1-x)^{\al+n}(1+x)^{\be+n}\right).
\]
\paragraph{Special cases}
\begin{itemize}
\item
{\em Gegenbauer} or {\em ultraspherical polynomials} ($\al=\be=\la-\thalf$):
\quad$\dstyle C_n^\la(x):=
\frac{(2\la)_n}{(\la+\thalf)_n}\,P_n^{(\la-\half,\la-\half)}(x)$.
\item
{\em Legendre polynomials} ($\al=\be=0$):\quad
$P_n(x):=P_n^{(0,0)}(x)$.
\item

{\em Chebyshev polynomials} ($\al=\be=\pm\thalf$):
\[
T_n(\cos\tha):=\cos(n\,\tha)=\frac{n!}{(\thalf)_n}\,P_n^{(-\half,-\half)}(\cos\tha),
\quad
U_n(\cos\tha):=\frac{\sin(n+1)\tha}{\sin\tha}=\frac{(2)_n}{(\tfrac32)_n}\,
P_n^{(\half,\half)}(\cos\tha).
\]
\end{itemize}
\paragraph{Quadratic transformations}
Since $P_{2n}^{(\al,\al)}(x)$ is an even polynomial of degree
$2n$ in $x$, it is also a polynomial $p_n(2x^2-1)$
of degree $n$ in $2x^2-1$. For $m\ne n$ we have
\[
0=\int_0^1 p_m(2y^2-1)p_n(2y^2-1)\,(1-y^2)^\al\,dy
=\const\int_{-1}^1p_m(x)p_n(x)\,(1-x)^\al(1+x)^{-\half}\,dx.
\]
Hence
\[
\frac{P_{2n}^{(\al,\al)}(x)}{P_{2n}^{(\al,\al)}(1)}=
\frac{P_n^{(\al,-\half)}(2x^2-1)}{P_n^{(\al,-\half)}(1)}.
\]
Similarly,
\[
\frac{P_{2n+1}^{(\al,\al)}(x)}{P_{2n+1}^{(\al,\al)}(1)}=
\frac{xP_n^{(\al,\half)}(2x^2-1)}{P_n^{(\al,\half)}(1)}\,.
\]
\begin{theorem}{\rm \cite[Ch.~1, \S8]{1}}\quad
Let $\{p_n\}$ be a system of
orthogonal polynomials with respect to an even weight function
$w$ on $\RR$. Then there are systems $\{q_n\}$ and $\{r_n\}$ of orthogonal
polynomials on $[0,\iy)$ with respect
to weight functions
$x\mapsto x^{-\half}w(x^\half)$ and $x\mapsto x^{\half}w(x^\half)$,
respectively, such that
$p_{2n}(x)=q_n(x^2)$ and $p_{2n+1}(x)=x\,r_n(x^2)$.\end{theorem}

\subsection{Electrostatic interpretation of zeros}
Let $p_n(x):=\const P_n^{(2p-1,2q-1)}(x)=(x-x_1)(x-x_2)\ldots(x-x_n)$ be monic
Jacobi polynomials ($p,q>0$). By \eqref{30}
\begin{equation*}
(1-x^2)p_n''(x)+2(q-p-(p+q)x)p_n'(x)=-n(n+2p+2q-1)p_n(x).
\end{equation*}
Hence
\begin{align*}
0&=1-x_k^2)p_n''(x_k)+2(q-p-(p+q)x_k)p_n'(x_k)\\
&=\thalf\,\frac{p_n''(x_k)}{p_n'(x_k)}+\frac p{x_k-1}+\frac q{x_k+1}
=\sum_{j,\;j\ne k}\frac1{x_k-x_j}+\frac p{x_k-1}+\frac q{x_k+1}.
\end{align*}
This can be reformulated as
\begin{equation*}
(\nabla V)(x_1,\ldots,x_n)=0,
\end{equation*}
where
\begin{equation*}
V(y_1,\ldots,y_n)
:=-\sum_{i<j}\log(y_j-y_i)-p\sum_j\log(1-y_j)-q\sum_j\log(1+y_j)
\end{equation*}
is the logarithmic potential obtained from $n+2$ charges $q,1,\ldots,1,p$ at
successive points
$-1<y_1<\allowbreak\ldots<y_n<1$. It
achieves a minimum at the zeros of $P_n^{(2p-1,2q-1)}$.
This result goes back to Stieltjes \cite{12}.

\subsection{Laguerre polynomials}
{\em Laguerre polynomials} $L_n^{\raisemath{1pt}\al}$
\cite[\S9.12]{4}
are orthogonal on $[0,\iy)$ with respect to the weight function
$w(x):=x^\al e^{-x}$ ($\al>-1$). They are normalized by
$L_n^{\raisemath{1pt}\al}(0)=(\al+1)_n/n!\,$.
They can be expressed in terms of terminating confluent hypergeometric
functions by
\begin{equation}
L_n^{\raisemath{1pt}\al}(x)=\frac{(\al+1)_n}{n!}\,
\hyp11{-n}{\al+1}x
=\sum_{k=0}^n\frac{(\al+k+1)_{n-k}}{k!\,(n-k)!}\,(-x)^k.
\label{32}
\end{equation}

For $p_n(x):=L_n^{\raisemath{1pt}\al}(x)$
there is the second order differential equation
\[
x\,p_n''(x)+(\al+1-x)\,p_n'(x)=-n\,p_n(x).
\]
This can be split up by the shift operator relations
\[
\frac d{dx}\,L_n^{\raisemath{1pt}\al}(x)=-L_{n-1}^{\raisemath{1pt}{\al+1}}(x),
\]
and
\begin{equation}
x\,\frac d{dx}\,L_{n-1}^{\raisemath{1pt}{\al+1}}(x)
+(\al+1-x)L_{n-1}^{\raisemath{1pt}{\al+1}}(x)=
x^{-\al}e^x\,\frac d{dx}\left(x^{\al+1}e^{-x}
L_{n-1}^{\raisemath{1pt}{\al+1}}(x)\right)=n\,L_n^{\raisemath{1pt}\al}(x).
\label{3}
\end{equation}
Iteration of \eqref{3} gives the Rodrigues formula
\[
L_n^{\raisemath{1pt}\al}(x)=\frac{x^{-\al}e^x}{n!}\,
\left(\frac d{dx}\right)^n\left(x^{n+\al}e^{-x}\right).
\]

\subsection{Hermite polynomials}
Hermite polynomials $H_n$ \cite[\S9.15]{4}
are orthogonal with respect to the weight function
$w(x):=e^{-x^2}$ on $\RR$ and they are normalized by
$H_n=2^n x^n+\cdots\;$.
They have the explicit expression
\begin{equation}
H_n(x)=n!\,\sum_{j=0}^{[n/2]}\frac{(-1)^j(2x)^{n-2j}}{j!\,(n-2j)!}\,.
\label{24}
\end{equation}
There is the second order differential equation
\[
H_n''(x)-2xH_n'(x)=-2nH_n(x).
\]
This can be split up by the shift operator relations
\begin{equation}
H_n'(x)=2n\,H_{n-1}(x),\quad
H_{n-1}'(x)-2xH_{n-1}(x)
=e^{x^2}\,\frac d{dx}\left(e^{-x^2}H_{n-1}(x)\right)
=-H_n(x).
\label{4}
\end{equation}
Iteration of the last equality in \eqref{4} gives the Rodrigues formula
\[
H_n(x)=(-1)^n\,e^{x^2}\left(\frac d{dx}\right)^n\left(e^{-x^2}\right).
\]
\subsection{General method to derive the standard formulas}
The previous formulas can be derived by the following general method.
Let $(a,b)$ be an open interval and let 
$w,w_1>0$ be strictly positive $C^1$-functions on $(a,b)$.
Let $\{p_n\}$ and $\{q_n\}$ be systems of
monic orthogonal polynomials on $(a,b)$
with respect to the weight function $w$ resp.\ $w_1$.
Then under suitable boundary assumptions for $w$ and $w_1\,$ we have
\begin{equation*}
\int_a^b p_n'(x)\,q_{m-1}(x)\,w_1(x)\,dx
=-\int_a^b p_n(x)\,w(x)^{-1}\,\frac d{dx}\big(w_1(x)\,q_{m-1}(x)\big)\,w(x)\,dx.
\end{equation*}
Suppose that for certain $a_n\ne0\,$:
\[
w(x)^{-1}\,\frac d{dx}\left(w_1(x)\,x^{n-1}\right)=-a_n\,x^n+
\mbox{polynomial of degree $<n$.}
\]
Then we easily derive a pair of first order differentiation formulas connecting
$\{p_n\}$ and $\{q_n\}$, an eigenvalue equation for $p_n$ involving a second
order differential operator, and a formula connecting the quadratic norms for
$p_n$ and $q_{n-1}\,$:
\begin{align*}
&p_n'(x)=n\,q_{n-1}(x),\quad
w(x)^{-1}\,\frac d{dx}\left(w_1(x)\,q_{n-1}(x)\right)=-a_n\,p_n(x),
\\
&w(x)^{-1}\,\frac d{dx}\left(w_1(x)\,p_n'(x)\right)=-na_n\,p_n(x),
\\
&n\int_a^b q_{n-1}(x)^2\,w_1(x)\,dx=a_n\int_a^b p_n(x)^2\,w(x)\,dx.
\end{align*}

In particular, if we work
with monic Jacobi polynomials
$p_n^{(\al,\be)}$, then
$(a,b)=(-1,1)$,
$w(x)=(1-x)^\al(1+x)^\be$, $p_n(x)=p_n^{(\al,\be)}(x)$,
$w_1(x)=(1-x)^{\al+1}(1+x)^{\be+1}$,
$q_m(x)=p_m^{(\al+1,\be+1)}(x)$.
Then $a_n=(n+\al+\be+1)$. Hence
\begin{align}
&\frac d{dx}\,p_n^{(\al,\be)}(x)=n\,p_{n-1}^{(\al+1,\be+1)}(x),
\label{12}\\
&\Big((1-x^2)\,\frac d{dx}
+\big(\be-\al-(\al+\be+2)x\big)\Big)p_{n-1}^{(\al+1,\be+1)}(x)
=-(n+\al+\be+1)\,p_n^{(\al,\be)}(x).
\label{13}
\end{align}
For $x=1$ \eqref{13} yields
\[
p_n^{(\al,\be)}(1)=\frac{2(\al+1)}{n+\al+\be+1}\,p_{n-1}^{(\al+1,\be+1)}(1).
\]
Iteration gives
\begin{equation}
p_n^{(\al,\be)}(1)=\frac{2^n(\al+1)_n}{(n+\al+\be+1)_n}.
\label{14}
\end{equation}
So for $p_n=\const p_n^{(\al,\be)}=k_n x^n+\cdots\;$
we know $p_n(1)/k_n$, independent of the normalization.

The hypergeometric series representation of Jacobi polynomials
is next obtained from \eqref{12} by Taylor expansion:
\begin{multline*}
\frac{p_n^{(\al,\be)}(x)}{p_n^{(\al,\be)}(1)}=
\sum_{k=0}^n\frac{(x-1)^k}{k!}\,\left(\frac d{dx}\right)^k p_n^{(\al,\be)}(x)\Big|_{x=1}\\
=
\sum_{k=0}^n\frac{(x-1)^k}{k!}\,\frac{n!}{(n-k)!}\,
\frac{p_{n-k}^{(\al+k,\be+k)}(1)}{p_n^{(\al,\be)}(1)}
=\hyp21{-n,n+\al+\be+1}{\al+1}{\thalf(1-x)}.
\end{multline*}

The quadratic norm $h_n$ can be obtained by iteration of
\begin{equation*}
\int_{-1}^1 p_n^{(\al,\be)}(x)^2\,(1-x)^\al(1+x)^\be\,dx
=\frac n{n+\al+\be+1}\int_{-1}^1
p_{n-1}^{(\al+1,\be+1)}(x)^2\,(1-x)^{\al+1}(1+x)^{\be+1}\,dx.
\end{equation*}
So for $p_n=\const p_n^{(\al,\be)}=k_n x^n+\cdots\;$
we know $h_n/k_n^2$, independent of the normalization.

\subsection{Characterization theorems}
Up to a constant factors and up to
transformations $x\to ax+b$ of the argument
the very classical orthogonal polynomials (Jacobi, Laguerre and Hermite)
are uniquely
determined as orthogonal polynomials
$p_n$ satisfying any of the following three criteria.
(In fact there are more ways to characterize these polynomials,
see Al-Salam \cite{13}.)
\begin{itemize}
\item ({\em Bochner theorem})\quad
The $p_n$ are eigenfunctions of a second order differential operator.
\item
The polynomlals $p_{n+1}'$ are again orthogonal polynomials.
\item
The polynomials are orthogonal with respect to a positive $C^\iy$ weight
function $w$ on an open interval $I$ and there is a polynomial $X$
such that the {\em Rodrigues formula} holds on~$I$:
\[
p_n(x)=\const w(x)^{-1}\,\left(\frac d{dx}\right)^n \big(X(x)^n w(x)\big).
\]
\end{itemize}

\subsection{Limit results}
The very classical orthogonal polynomials are connected to each other by
limit relations. We give these limits below for the monic versions
$p_n^{(\al,\be)}, \ell_n^\al, h_n$ of these polynomials, and on each line
we give also the corresponding limit for the weight functions:
\begin{align}
\al^{n/2}p_n^{(\al,\al)}(x/\al^{1/2})&\to h_n(x),
&(1-x^2/\al)^\al\to e^{-x^2},\quad\al\to\iy,\label{26}\\
(-\be/2)^n\,p_n^{(\al,\be)}(1-2x/\be)&\to\ell_n^\al(x),
&x^\al(1-x/\be)^\be\to x^\al e^{-x},\quad\be\to\iy,\label{27}\\
(2\al)^{-n/2}\,\ell_n^\al((2\al)^{1/2} x+\al)&\to h_n(x),&
(1+(2/\al)^{1/2}x)^\al e^{-(2\al)^{1/2} x}\to e^{-x^2},\quad\al\to\iy.\label{28}
\end{align}
The limits of the orthogonal polynomials in \eqref{26} and \eqref{27}
immediately follow from \eqref{29}, \eqref{32} and \eqref{24}.
For various ways to prove \eqref{28} see
\cite[section 2]{14}.
\section{Gauss quadrature}
Let be given $n$ real points $x_1<x_2<\ldots<x_n$.
Put $p_n(x):=(x-x_1)\ldots(x-x_n)$. For $k=1,\ldots,n$
let $l_k$ be the unique polynomial of degree $<n$ such that
$l_k(x_j)=\de_{k,j}\,$ ($j=1,\ldots,n$). This polynomial,
called the {\em Lagrange interpolation polynomial}, equals
\[
l_k(x)=\frac{\prod_{j;\,j\ne k}(x-x_j)}{\prod_{j;\,j\ne k}(x_k-x_j)}
=\frac{p_n(x)}{(x-x_k)\,p_n'(x_k)}\,.
\]
For all polynomials $r$ of degree $<n$ we have
\[
r(x)=\sum_{k=1}^n r(x_k)\,l_k(x).
\]
\begin{theorem}[Gauss quadrature]
Let $p_n$ be an orthogonal polynomial with respect to $\mu$
and let the $l_k$
be the Lagrange interpolation polynomials associated with the zeros
$x_1,\ldots,x_n$ of $p_n$.. Put
\[
\la_k:=\int_\RR l_k(x)\,d\mu(x).
\]
Then
\[
=\int_\RR l_k(x)^2\,d\mu(x)>0
\]
and for all polynomials of degree $\le 2n-1$ we have
\begin{equation}
\int_\RR f(x)\,d\mu(x)=\sum_{k=1}^n \la_k\,f(x_k).
\label{8}
\end{equation}
\end{theorem}
\Proof
Let $f$ be a polynomial of degree $\le 2n-1$.
Then for certain polynomials $q$ and $r$ of degree $\le n-1$ we have
$f(x)=q(x)p_n(x)+r(x)$. Hence $f(x_k)=r(x_k)$ and
\begin{equation*}
\int_\RR f(x)\,d\mu(z)=\int_\RR r(x)\,d\mu(x)
=\sum_{k=1}^n r(x_k)\,\int_\RR l_k(x)\,d\mu(x)
=\sum_{k=1}^n \la_k r(x_k)=\sum_{k=1}^n\la_k f(x_k).
\end{equation*}
Also
\[
\qquad\qquad\qquad\qquad\qquad
\la_k=\sum_{j=1}^n\la_j\,l_k(x_j)^2=\int_\RR l_k(x)^2\,d\mu(x)>0.
\qquad\qquad\qquad\qquad\qquad\qquad\qed
\]

From \eqref{8} we see in particular that, for $i,j\le n-1$,
\[
h_j\de_{i,j}=\int_\RR p_i(x)p_j(x)\,d\mu(x)=
\sum_{k=1}^n\la_k\,p_i(x_k)\,p_j(x_k).
\]
Thus the finite system $p_0,p_1,\ldots,p_{n-1}$ forms a set of orthogonal
polynomials on the finite set $\{x_1,\ldots,x_n\}$
of the $n$ zeros of $p_n$ with respect to
the weights $\la_k$ and with quadaratic norms~$h_j$.
All information about this system is already contained in the finite system of
recurrence relations
\[
xp_j(x)=a_jp_{j+1}(x)+b_jp_j(x)+c_jp_{j-1}(x)\quad(j=0,1,\ldots,n-1)
\]
with $a_jc_{j+1}>0$ ($j=0,1,\ldots,n-2$). In particular, the $\la_k$ are
obtained up to a constant factor by solving the system
\[
\sum_{k=1}^n\la_k p_j(x_k)=0\quad(j=1,\ldots,n-1).
\]

\section{Askey scheme}
As an example of a finite system of orthogonal polynomials as described at the
end of the previous section,
consider orthogonal polynomials $p_0,p_1,\ldots,p_N$
on the zeros $0,1,\ldots,N$ of the polynomial
$p_{N+1}(x):=x(x-1)\ldots(x-N)$ with respect to nice explicit
weights $w_x$ ($x=0,1,\ldots,N$) like:
\begin{itemize}
\item
$w_x:=\binom nx p^x(1-p)^{N-x}$\quad ($0<p<1$).
Then the $p_n$ are the {\em Krawtchouk polynomials}
\[
K_n(x;p,N):=\hyp21{-n,-x}{-N}{\frac 1p}
=\sum_{k=0}^n\frac{(-n)_k(-x)_k}{(-N)_k\,k!}\,\frac1{p^k}\,.
\]
\item
$\dstyle w_x:=\frac{(\al+1)_x}{x!}\,\frac{(\be+1)_{N-x}}{(N-x)!}$\quad
($\al,\be>-1$).
Then the $p_n$ are the {\em Hahn polynomials}
\[
Q_n(x;\al,\be,N):=\hyp32{-n,n+\al+\be+1,-x}{\al+1,-N}1.
\]
\end{itemize}

Hahn polynomials are discrete versions of Jacobi polynomials:
\begin{multline*}
Q_n(Nx;\al,\be,N)=\hyp32{-n,n+\al+\be+1,-Nx}{\al+1,-N}1\to\\
\hyp21{-n,n+\al+\be+1}{\al+1}x=\const P_n^{(\al,\be)}(1-2x)
\end{multline*}
and
\begin{multline*}
N^{-1}
\sum_{x\in\{0,\frac1N,\frac2N,\ldots,1\}}Q_m(Nx;\al,\be,N)Q_n(Nx;\al,\be,N)\,w_{Nx}\to\\
\const \int_0^1 P_m^{(\al,\be)}(1-2x)P_n^{(\al,\be)}(1-2x)\,
x^\al(1-x)^\be\,dx.
\end{multline*}

Jacobi and Krawtchouk polynomials are different ways of looking at
the matrix elements of the irreducible representations of ${\rm SU}(2)$, see \cite{16}.
The $3j$ coefficients or Clebsch-Gordan coefficients for ${\rm SU}(2)$
can be expressed as Hahn polynomials,
see for instance \cite{15}.

While we saw that the Jacobi, Laguerre and Hermite polynomials are
eigenfunctions of a second order differential operator,
\begin{equation}
A(x)p_n''(x)+B(x)p_n'(x)+C(x)p_n(x)=\la_n p_n(x),
\label{15}
\end{equation}
the Hahn and Krawtchouk polynomials are examples of
orthogonal polynomials $p_n$
on $\{0,1,\ldots,N\}$ which are eigenfunctions of a second order
difference operator,
\begin{equation}
A(x)p_n(x-1)+B(x)p_n(x)+C(x)p_n(x+1)=\la_n\,p_n(x).
\label{16}
\end{equation}

If we also allow orthogonal polynomials on the infinite set
$\{0,1,2,\ldots\}$ then
{\em Meixner polynomials} $M_n(x;\be,c)$ and
{\em Charlier polynomials} $C_n(x;a)$ appear. Here
\begin{align*}
M_n(x;\be,c)&:=\hyp21{-n,-x}\be{1-\,\frac1c},\quad
w_x:=\frac{(\be_x)}{x!}\,c^x,\\
C_n(x;a)&:={}_2F_0(-n,-x;;-a^{-1}),\quad w_x:=a^x/x!\,.
\end{align*}
If we also include orthogonal polynomials which are
eigenfunctions of a second order operator as follows,
\begin{equation}
A(x) p_n(x+i)+B(x)p_n(x)+C(x)p_n(x-i)=\la_n p_n(x),
\label{17}
\end{equation}
then we have collected all families of orthogonal polynomials which
belong to the {\em Hahn class}.

Similarly, with an eigenvalue equation of the form
\begin{equation}
A(x) p_n(q(x+1))+B(x) p_n(q(x))+C(x) p_n(q(x-1))=\la_n p_n(q(x)),
\label{18}
\end{equation}
where $q$ is a fixed polynomial of second degree, we obtain
the orthogonal polynomials on  a {\em quadratic lattice}.
All orthogonal polynomials satisfying an equation of the form \eqref{15}--\eqref{18}
have been classified.
There are only 13 families, depending on at most four parameters,
and all expressible as hypergeometric functions,
${}_4F_3$ in the most complicated case.
They can be arranged hierarchically according to limit transitions denoted
by arrows. This is the famous {\em Askey scheme}, see for instance
\cite[Fig.1]{14}.

\section{The $q$-case}
On top of the Askey-scheme is lying the $q$-Askey scheme
\cite[beginning of Ch.~14]{4}, from which
there are also arrows to the Askey scheme as $q\to 1$.
We take always $0<q<1$ and let $q\uparrow 1$ for the limit to the classical case.
Some typical examples of $q$-analogues of classical concepts are
(see Gasper \& Rahman \cite{17}):
\begin{itemize}
\item
$q$-number: $\dstyle [a]_q:=\frac{1-q^a}{1-q}\to a$
\item
$q$-shifted factorial: $\dstyle(a;q)_n:=\prod_{k=0}^{n-1}(1-aq^k)$
(also for $n=\iy$),\quad
$\dstyle\frac{(q^a;q)_k}{(1-q)^a}\to (a)_k$.
\item
$q$-hypergeometric series:
\[
\qhyp{s+1}s{a_1,\ldots,a_{s+1}}{b_1,\ldots,b_s}{q,z}:=
\sum_{k=0}^\iy \frac{(a_1;q)_k\ldots(a_{s+1};q)_k}{(b_1;q)_k\ldots(b_s;q)_k}\,
\frac{z^k}{(q;q)_k}\,,
\]
\[
\qhyp{s+1}s{q^{a_1},\ldots,q^{a_{s+1}}}{q^{b_1},\ldots,q^{b_s}}{q,z}\to
\hyp{s+1}s{a_1,\ldots,a_{s+1}}{b_1,\ldots,b_s}z.
\]
\item
$q$-derivative:
$\dstyle (D_qf)(x):=\frac{f(x)-f(qx)}{(1-q)x}\to f'(x)$.
\item
$q$-integral:
$\dstyle \int_0^1 f(x)\,d_qx:=(1-q)\sum_{k=0}^\iy f(q^k)\,q^k\to\int_0^1 f(x)\,
dx$.
\end{itemize}

The $q$-case allows more symmetry which may be broken when taking limits
for $q$ to 1. In the elliptic case \cite[Ch.~11]{17}
lying above the $q$-case there is even more
symmetry.

On the highest level in the $q$-Askey scheme are the
{\em Askey--Wilson polynomials} \cite{18}. They are given by
\begin{equation*}
p_n(\cos\tha;a,b,c,d\mid q):=\frac{(ab;q)_n(ac;q)_n(ad;q)_n}{a^n}\,
\qhyp43{q^{-n},q^{n-1}abcd,ae^{i\tha},ae^{-i\tha}}{ab,ac,ad}{q,q},
\end{equation*}
and they are symmetric in the parameters $a,b,c,d$.
For suitable restrictions on the parameters
they are orthogonal with respect to an explicit weight function on $(-1,1)$.
In the special case
$a=-c=\be^\half,\;b=-d=(q\be)^\half$
we get the {\em continuous $q$-ultraspherical polynomials} \cite[\S14.10.1]{4}.
They satisfy the orthogonality relation
\[
\int_0^\pi C_m(\cos\tha;\be\mid q)\,C_n(\cos\tha;\be\mid q)\,
\left|\frac{(e^{2i\tha};q)_\iy}{(\be e^{2i\tha};q)_\iy}\right|^2\,d\tha=0\qquad
(m\ne n),
\]
and they have the generating function
\[
\left|\frac{(\be e^{i\tha}t;q)_\iy}{(e^{i\tha}t;q)_\iy}\right|^2
=\sum_{n=0}^\iy C_n(x;\be\mid q)t^n.
\]
For $q\uparrow 1$ they tend to ultraspherical (or Gegenbauer) polynomials:
$C_n(x;q^\la\mid q)\to C_n^\la(x)$.
The Gegenbauer polynomials have the generating function
\[
(1-2xt+t^2)^{-\la}=\sum_{n=0}^\iy C_n^\la(x) t^n.
\]

\section{Some deeper properties of general orthogonal polynomials}
\subsection{True interval of orthogonality}
Consider a system of orthogonal polynomials $\{p_n\}$.
Let $p_n$ have zeros $x_{n,1}<x_{n,2}<\ldots<x_{n,n}\,$.
Then
\begin{equation}
\begin{split}
x_{i,i}>x_{i+1,i}>\ldots>x_{n,i}&\downarrow\xi_i\ge-\iy,\\
x_{j,1}<x_{j+1,2}<\ldots<x_{n,n-j+1}&\uparrow\eta_j\le\iy.
\end{split}
\label{33}
\end{equation}
The closure $I$ of the interval $(\xi_1,\eta_1)$ is called the
{\em true interval of orthogonality} of the system $\{p_n\}$
(see \cite[Ch.~I, Definition 5.2]{1}.
It has the following properties.
\begin{enumerate}
\item
$I$ is the smallest closed interval containing all zeros $x_{n,i}$.
\item
There is an orthogonality measure $\mu$ for the $p_n(x)$
such that $I$ is the
smallest closed interval containing the support of $\mu$
(see \cite[p.113]{7} or \cite[Ch.~II, Theorem 4.2]{1})
\item
If $\mu$ is any orthogonality measure for the $p_n(x)$ and $J$
is a closed interval containing the support of $\mu$ then $I\subset J$
(see \cite[Ch.~II, Theorem 4.1]{1}).
\end{enumerate}

\subsection{Criteria for bounded support of orthogonality measure}
Recall the three-term recurrence relation \eqref{6} for a system of
monic orthogonal polynomials~$\{p_n\}$.
Let $\xi_1,\eta_1$ be as in \eqref{33}.
The following theorem gives criteria
for the support of an orthogonality measure in terms of the behaviour of the
coefficients $b_n, c_n$ in \eqref{6} as $n\to\iy$.
\begin{theorem}\quad
\begin{enumerate}
\item
{\rm(\cite[p.109]{1})}
If $\{b_n\}$ is bounded and $\{c_n\}$ is unbounded then
$(\xi_1,\eta_1)=(-\iy,\iy)$.
\item
{\rm(\cite[p.109, Theorem 2.2]{1})}
$\{b_n\}$ and $\{c_n\}$ are bounded iff
$[\xi_1,\eta_1]$ is bounded.
\item
{\rm(\cite[Theorem 10]{26})}
If $b_n\to b$ and $c_n\to c$ {\rm (}$b,c$ finite{\rm )} then
{\rm supp}$(\mu)\supset[b-2\sqrt c,b+2\sqrt c\,]$, and outside
$[b-2\sqrt c,b+2\sqrt c\,]$ the support of $\mu$ has at most
countably many points $x_k$, while the only possible limit points of
$\{x_k\}$ are $b\pm 2\sqrt c$.

\end{enumerate}
\end{theorem}
\begin{example}
\rm
Monic Jacobi polynomials $p_n^{(\al,\be)}(x)\,$:
\begin{align*}
b_n&=\frac{\be^2-\al^2}
{(2n+\al+\be)(2n+\al+\be+2)}\,\to 0,\sLP
c_n&=\frac{4n(n+\al)(n+\be)(n+\al+\be)}
{(2n+\al+\be-1)(2n+\al+\be)^2(2n+\al+\be+1)}\,\to\frac14\,.
\end{align*}
Hence $[b-2\sqrt c,b+2\sqrt c\,]=[-1,1]$.
\end{example}

\subsection{Criteria for uniqueness of orthogonality measure}
Put
\[
\rho(z):=\Big(\sum_{n=0}^\iy |p_n(z)|^2\Big)^{-1}\quad(z\in\CC).
\]
Then $0\le\rho(z)<\iy$.  Note that $\rho(z)=0$ iff 
$\sum_{n=0}^\iy |p_n(z)|^2$ diverges and that $\rho(z)>0$ iff
$\sum_{n=0}^\iy |p_n(z)|^2$ converges.
\begin{theorem}
{\rm(\cite[pp.~49--51]{7})}
The orthogonality measure is not unique iff
$\rho(z)>0$ for all $z\in\CC$.
Equivalently, the orthogonality measure is unique iff $\rho(z)=0$ for some $z\in\CC$.
\sLP
In the case of a unique orthogonality measure $\mu$, we have $\rho(z)=0$
for $z\in\CC\backslash\RR$ and
$\rho(x)=\mu(\{x\})$ {\rm (}the mass at $x${\rm )} for $x\in\RR$,
which implies that $\rho(x)\ne0$ iff $x$ is a mass point of $\mu$.
\sLP
In case of non-uniqueness, for each $x\in\RR$ the largest possible mass
of a measure $\mu$ at $x$ is $\rho(x)$ and there is a measure realizing this mass
at $x$.
\end{theorem}

Recall the moments $\mu_n:=\lan x^n,1\ran=\int_\RR x^n\,d\mu(x)$,
which are uniquely determined (up to a constant factor) by the system $\{p_n\}$,
and also recall the three-term recurrence relation \eqref{7} for a system of
orthonormal polynomials $\{p_n\}$. 

\begin{theorem}[Carleman]
{\rm(\cite[Theorem 1.10 and pp.~47, 59]{7})}
\label{23}
There is a unique orthogonality measure for the $p_n$ if one of the following two
conditions is satisfied.
\mPP
{\rm(i)}\quad $\dstyle \sum_{n=1}^\iy \mu_{2n}^{-1/(2n)}=\iy$,\qquad
{\rm(ii)}\quad $\dstyle\sum_{n=1}^\iy a_n^{-1}=\iy$.
\end{theorem}

\begin{example}[Hermite]
\rm
\[
\mu_{2n}=\int_{-\iy}^\iy x^{2n}e^{-x^2}\,dx=\Ga(n+\thalf)
\quad
{\rm and}
\quad
\log\Gamma(n+\thalf)=n\log(n+\thalf)+O(n)\quad\mbox{as $n\to\iy$},
\]
so $\mu_{2n}^{-1/(2n)}\sim(n+\thalf)^{-\half}$.
Hence $\dstyle \sum_{n=1}^\iy \mu_{2n}^{-1/(2n)}=\iy\,$, i.e.,
the orthogonality measure for the Hermite polynomials is unique.
\end{example}

\begin{example}[Laguerre]
\rm
Monic Laguerre polynomials $p_n$ satisfy
\[
xp_n(x)=p_{n+1}(x)+(2n+\al+1)p_n(x)+n(n+\al)p_{n-1}(x).
\]
Since
$\dstyle\sum_{n=0}^\iy\frac1{(n(n+\al))^{1/2}}=\iy\,$, 
the orthogonality measure is unique.
Also note that
\[
\frac{L_n^{\raisemath{1pt}\al}(0)^2}{h_n}=\frac{\Ga(n+\al+1)}{\Ga(n+1)\Ga(\al+1)}
\sim n^\al.
\]
Since
$\dstyle \sum_{n=1}^\iy n^\al=\iy$\quad($\al>-1$) we conclude once more
that the orthogonality measure is unique.
\end{example}

\begin{example}[Stieltjes--Wigert]
\rm
Consider the moments $\mu_n$ given by the \RHS\ of \eqref{22}. Then
\[
\sum_{n=1}^\iy \mu_{2n}^{-1/(2n)}=
\sum_{n=1}^\iy e^{-\half(n+1)}<\iy.
\]
Since the
corresponding moment problem is undetermined, the above inequality agrees
with Theorem \ref{23}(i). Furthermore, from \cite[(14.27.4)]{4} with
$q=e^{-\half}$ we see that the corresponding orthonormal polynomials
$p_n(x)=\const S_n(q^\half x;q)$ with $q=e^{-\half}$
have $a_{n-1}^2=e^{2n}(1-e^{-\half n})$, by which
$\sum_{n=1}^\iy a_n^{-1}<\iy$, in agreement with Theorem \ref{23}(ii).

\end{example}

\subsection{Orthogonal polynomials and continued fractions}
Let monic orthogonal polynomials $p_n$ be recursively defined by
\[
p_0(x)=1,\qquad p_1(x)=x-b_0,\qquad
xp_n(x)=p_{n+1}(x)+b_np_n(x)+c_np_{n-1}(x)\quad(n\ge 1,\;c_n>0).
\]
Then define monic orthogonal polynomials $p_n^{(1)}$ by
\[
p_0^{(1)}(x)=1,\qquad p_1^{(1)}(x)=x-b_1,\qquad
xp_n^{(1)}(x)=p_{n+1}^{(1)}(x)+b_{n+1}p_n^{(1)}(x)+c_{n+1}p_{n-1}^{(1)}(x)
\quad(n\ge 1).
\]
They are called {\em first associated orthogonal polynomials} or
{\em numerator polynomials}.

Define
\[
F_1(x):=\frac1{x-b_0}\,,\quad
F_2(x):=\frac1{x-b_0-\frac{c_1}{x-b_1}}\,,\quad
F_3(x):=\frac1{x-b_0-\frac{c_1}{x-b_1-\frac{c_2}{x-b_2}}}\,,
\]
and recursively obtain $F_{n+1}(x)$ from $F_n(x)$ by replacing $b_{n-1}$ by
$\dstyle b_{n-1}+\frac{c_n}{x-b_n}\,$.
This is a {\em continued fraction}, which can be notated as
\[
F_n(z)=\frac1{z-b_0-|}\;\frac{|c_1}{z-b_1-|}\;\cdots\,\frac{|c_{n-2}}{z-b_{n-2}-|}\;
\frac{|c_{n-1}}{z-b_{n-1}}\,.
\]
\begin{theorem}[essentially due to Stieltjes]
{\rm(\cite[Ch.~3, \S4]{1})}
\[
F_n(z)=\frac{p_{n-1}^{(1)}(z)}{p_n(z)}\qquad{\rm and}\qquad
p_{n-1}^{(1)}(y)=\frac1{\mu_0}\,
\int_\RR\frac{p_n(y)-p_n(x)}{y-x}\,d\mu(x).
\]
\end{theorem}

\begin{theorem}[Markov]
{\rm(\cite[Ch.~3, (4.8)]{1})}
Suppose that there is a {\rm(}unique{\rm)} orthogonality measure $\mu$ of
bounded support
for the $p_n$.
Let $[\xi_1,\eta_1]$ be the true interval of orthogonality. Then
\[
\lim_{n\to\iy} F_n(z)=\frac1{\mu_0}\int_{\xi_1}^{\eta_1}\frac{d\mu(x)}{z-x},
\]
uniformly on compact subsets of $\CC\backslash[\xi_1,\eta_1]$.
\end{theorem}

\subsection{Measures in case of non-uniqueness}
Take $p_n$ and $p_n^{(1)}$ orthonormal:
\begin{align*}
&p_0(x)=1,\quad p_1(x)=(x-b_0)/a_0,\\
&xp_n(x)=a_np_{n+1}(x)+b_np_n(x)+a_{n-1}p_{n-1}(x)\quad(n\ge 1),\sLP
&p_0^{(1)}(x)=1,\quad p_1^{(1)}(x)=(x-b_1)/a_1,\\
&xp_n^{(1)}(x)=a_{n+1}p_{n+1}^{(1)}(x)+b_{n+1}p_n^{(1)}(x)+a_np_{n-1}^{(1)}(x)
\quad(n\ge 1),
\end{align*}
where $a_n>0$.
Let $\mu_0=1,\;\mu_1,\mu_2,\ldots\;$ be the moments for the $p_n$.
Suppose that the orthogonality measure for the $p_n$ is not unique.
Then the possible orthogonality measures are precisely the positive measures
$\mu$ solving the moment problem \eqref{19}.
The set of these solutions is convex and weakly compact.

We will need the following entire analytic functions.
\begin{align*}
A(z)&:=z\sum_{n=0}^\iy p_n^{(1)}(0)\,p_n^{(1)}(z),\quad
B(z):=-1+z\sum_{n=1}^\iy p_{n-1}^{(1)}(0)\,p_n(z),\\
C(z)&:=1+z\sum_{n=1}^\iy p_n(0)\,p_{n-1}^{(1)}(z),\quad
D(z)=z\sum_{n=0}^\iy p_n(0)\,p_n(z).
\end{align*}
By a {\em Pick function} we mean a holomorphic function $\phi$
mapping the open upper half plane into the closed upper half plane.
Let $\bf P$ denote the set of all Pick functions.
In the theorem below we will associate with a Pick function $\phi$ a certain
measure $\mu_\phi$. There $\mu_t$ for $t\in\RR$ will mean the measure
$\mu_\phi$ with $\phi$ the constant Pick function $z\mapsto t$, and
$\mu_\iy$ will mean the measure $\mu_\phi$ with $\phi$ the constant 
function $z\mapsto \iy$ (not a Pick function).

\begin{theorem}[Nevanlinna, M. Riesz]
{\rm(\cite[Theorem 2.12]{7})}
Suppose the moment problem \eqref{19} is undetermined.
The identity
\[
\int_\RR \frac{d\mu_\phi(t)}{t-z}=
-\,\frac{A(z)\phi(z)-C(z)}{B(z)\phi(z)-D(z)}\quad(\Im z>0)
\]
gives a one-to-one correspondence $\phi\to\mu_\phi$
between ${\bf P}\union\{\iy\}$ and
the set of measures solving the moment problem \eqref{19}.

Furthermore the measures $\mu_t$ $(t\in\RR\union\{\iy\})$ are precisely
the extremal elements of the convex set, and also precisely the measures $\mu$
solving \eqref{19} for which the the polynomials are dense in $L^2(\mu)$.
All measures $\mu_t$ are discrete. The mass points of $\mu_t$
are the zeros of the entire function $tB-D$ {\rm(}or of $B$ if $t=\iy${\rm)}.
\end{theorem}

\begin{example}[Stieltjes--Wigert]
\rm
The measure which gives in \eqref{25} a solution for the moment
problem associated with special Stieltjes--Wigert polynomials, has a support which
is almost discrete, but not completely, since $0$ is a limit point of the support.
Therefore (see the above theorem)
this measure cannot be extremal. As observed by Christiansen
at the end of \cite{9}, finding explicit extremal measures for this case seems
to be completely out of reach. Since the measure in \eqref{25} is not extremal,
the polynomials will not be dense in the corresponding $L^2$ space.
Christiansen \& Koelink \cite[Theorem 3.5]{10} give an explicit orthogonal
system in this $L^2$ space which complements the orthogonal system of
Stieltjes--Wigert polynomials to a complete orthogonal system.
\end{example}

\section{Orthogonal polynomials in connection with computer algebra}
Undoubtedly, computer algebra is nowadays a powerful tool which many
mathematicians and physicists use in daily practice for their research,
often using wide spectrum computer algebra programs
like Mathematica or Maple, to which further specialized packages
are possibly added. This is certainly also
the case for research in orthogonal polynomials, in particular when it concerns
special families. Jacobi, Laguerre and Hermite polynomials can be
immediately called in Mathematica and Maple, while other polynomials in
the ($q$-)Askey scheme can be defined by their
\mbox{($q$-)hypergeometric} series
interpretation. Even more geneneral special orthogonal polynomials
can be generated by their three-term recurrence relation.

Typical kinds of computations being done are:
\begin{enumerate}
\item
Checking a symbolic computation on computer which was first done by hand.
\item
Doing a symbolic computation first on computer and then find a hopefully elegant
derivation which can be written up.
\item
Doing a symbolic computation on computer and then write in the paper
something like: ``By using Mathematica we found $\ldots$''.
\item
Checking general theorems, with (hopefully correct) proofs available,
for special examples by computer algebra.
\item
Formulating general conjectures in interaction with output of symbolic
computation for special examples.
\item
Trying to find a simple evaluation of a parameter dependent expression
by extrapolating from outputs for special cases of the parameters.
\item
Building large collections of formulas, to be made available on the internet,
which are fully derived by computer algebra, and which can be made
adaptive for the user.
\item
Applying full force computer algebra, often using special purpose programs,
for obtaining massive output which is a priori hopeless to get by hand or to be
rewritten into an elegant expression.
\end{enumerate}

While item 8 is common practice in high energy physics, I have little to
say about this from my own experience. Concerning item 3 there may be a
danger that we become lazy, and no longer look for an elegant analytic proof
when the result was already obtained by computer algebra.
In particular, many formulas for terminating hypergeometric series can be derived
much quicker when we recognize them as orthogonal polynomials and use some
orthogonality argument.

As an example ot item 1, part of the formulas in the NIST handbook \cite{5}
was indeed checked by computer algebra.
Concerning item 7, it is certainly a challenge for computer algebrists how much of a formula database for special functions can be produced purely by
computer algebra. Current examples are CAOP \cite{19} (maintained by
Wolfram Koepf, Kassel) and DDMF \cite{20} (maintained by Fr\'ed\'eric Chyzak et al.\ at INRIA).

The most spectacular success of computer algebra for special functions has
been the Zeilberger algorithm, now already more than 20 years old.
It is treated in several books:
Petkov{\v s}ek et al.\ \cite{21}, Koepf \cite{22},
Kauers \& Paule \cite{23}. In particular, \cite{22} contains quite a lot
of examples of application of this algorithm to special orthogonal polynomials,
including the discrete and the $q$-case.

Various applications of computer algebra to special orthogonal polynomials
can be found in other chapters of the present volume.

A very desirable application of computer algebra would be to recognize
from a given three-term recurrence relation with explict, possibly still
parameter dependent coefficients, whether it comes from a system of
orthogonal polynomials in the \mbox{($q$-)Askey scheme}, and if so, which
system precisely. A very heuristic algorithm was implemented in
the procedure Rec2ortho \cite{24} (started by Swarttouw and maintained by
the author). It is only up to the level of 2 parameters in the Askey scheme.
On the other hand Koepf \& Schmersau \cite{25} give an algortithm
how to go back and forth between an explicit eigenvalue equation
\eqref{15} or \eqref{16} and a corresponding
three-term recurrence relation with
explicit coefficients.

\quad\\
\begin{footnotesize}
\begin{quote}
{ T. H. Koornwinder, Korteweg-de Vries Institute, University of
 Amsterdam,\\
 P.O.\ Box 94248, 1090 GE Amsterdam, The Netherlands;

\vspace{\smallskipamount}
email: }{\tt thkmath@xs4all.nl}
\end{quote}
\end{footnotesize}

\end{document}